\documentclass{amsart}
\usepackage[centertags]{amsmath}
\usepackage{graphicx}
\usepackage{amsfonts,amssymb,latexsym,epsfig,amsthm}
\usepackage{fullpage}

\newtheorem{lemma}{Lemma}

\newcommand{\ind}{1\!  {\rm l}}

\author{Jinho Baik  \and  Alexei Borodin \and  Percy Deift  \and  Toufic Suidan}
\address{Jinho Baik (baik@umich.edu): Department of Mathematics, University of Michigan,  Ann Arbor, MI}
\address{Alexei Borodin (borodin@caltech.edu): Department of Mathematics, California Institute of  Technology, Pasadena, CA}
\address{Percy Deift (deift@cims.nyu.edu): Courant Institute of Mathematical Sciences, New York  University, New York City, NY}
\address{Toufic Suidan (tsuidan@ucsc.edu): Department of Mathematics, University of California, Santa Cruz, CA}

\title{A Model for the Bus System in Cuernevaca (Mexico)}

\begin{document}

\maketitle

\section{Introduction}

The bus transportation system in Cuernevaca, Mexico, has certain distinguished, innovative features and has been the subject of an intriguing, recent study by M. \v{K}rbalek and P. \v{S}eba in~\cite{KS}.  The situation is as follows. We quote from~\cite{KS}: In Cuernevaca there...

\begin{center}
\parbox{5.5in}{
is no covering company responsible for organizing the city transport. Consequently, constraints such as a time table that represents external influence on the transport do not exist. Moreover, each bus is the property of the driver. The drivers try to maximize their income and hence the number of passengers they transport. This leads to competition among the drivers and to their mutual interaction. It is known that without additive interaction the probability distribution of the distances between subsequent buses is close to the Poisonian distribution and can be desribed by the standard bus route model\hspace{.02in}.\hspace{.02in}.\hspace{.02in}.\hspace{.05in}. A Poisson-like distribution implies, however, that the probability of close encounters of two buses is high (bus clustering) which is in conflict with the effort of the drivers to maximize the number of transported passengers and accordingly to maximize the distance to the preceding bus. In order to avoid the unpleasant clustering effect the bus drivers in Cuernevaca engage people who record the arrival times of buses at significant places.  Arriving at a checkpoint, the driver receives the information of when the previous bus passed that place. Knowing the time interval to the preceding bus the driver tries to optimize the distance to it by either slowing down or speeding up. In such a way the obtained information leads to a direct interaction between buses...}
\end{center}

In~\cite{KS}, \v{K}rbalek and \v{S}eba describe their work in analyzing the statistics of bus arrivals on Line 4 close to the city center. They study, in particular, the bus spacing distribution and also the bus number variance measuring the fluctuations of the total number of buses arriving at a fixed location during a time interval $T$. Quite remarkably, \v{K}rbalek and \v{S}eba find that these two statistics are well modeled by the Gaussian Unitary Ensemble (GUE) of random matrix theory (RMT)(see figures 2 and 3 in~\cite{KS}). Our goal in this paper is to provide a plausible explanation of these observations, and to this end we introduce a microscopic model for the bus line that leads simply and directly to GUE.


As noted in~\cite{KS} the number variance for the buses is in good agreement with the GUE formula up to a time interval T=3 (see figure 3 in \cite{KS}). As explained in~\cite{KS}, this behavior is consonant with the fact that each bus driver, using the information given by the recorder, interacts with the bus immediately behind him and the bus immediately in front of him. In other words, the primary interaction is a three body interaction. As is well known, particle systems modeled by GUE involve interactions between all the particles. This means, in particular, that if we use GUE to model a system with nearest neighbor (or more generally, short range) interactions we should restrict our attention to statistics that involve only nearest neighbor (or short range) interactions. This is the case for the spacing distributions and hence one is able to account for the good agreement between GUE and the observed data for this statistic across the entire parameter range in~\cite{KS}; this is in contrast to the number variance where, as noted above, there is good agreement only up to T=3. 


The paper is organized as follows. In Section \ref{Section2} we introduce our model. In Section \ref{Section3} we describe the double scaling limits of interest and indicate how to analyze these limits for our model using standard asymptotic techniques from random matrix theory. Our main results on the spacing distributions, as well as on the number variance, are stated at the end of Section \ref{Section3} . In Section \ref{Section4} we show how to modify the model to include alternative bus schedules.

Our analysis is based on the fact that certain nonintersecting paths
models lead to random matrix type ensembles, more exactly, orthogonal
polynomial ensembles. Various manifestations of this phenomenon can be
observed in~\cite{Johansson, KonigOConnell, KonigOConnellRoch, BorodinOlshanski}.

  Finally, in view of the calculations that follow, we note that formula
(4) in~\cite{KS} for the density of the spacing distribution is the so-called
"Wigner surmise" rather than the Gaudin distribution. The Gaudin
distribution (see  \eqref{conditional}  below) gives the exact formula for the spacing
distribution: However, the Wigner surmise is known to approximate this
exact formula to high accuracy (see Mehta~\cite{Mehta}). \vspace{.1in}

\noindent {\it Acknowledgments: The authors would like to thank Y. Avron for bringing the work of \v{K}rbalek and \v{S}eba to their attention. The authors would also like to thank P. \v{S}eba for providing useful information. The work of the authors was supported in part by NSF grants DMS0457335, DMS0402047, DMS0296084/0500923, and DMS0202530/0553403, respectively. }

\section{Basic Model}\label{Section2}

The dynamics of our bus model takes place on the lattice  $\mathbb{Z}$ as follows. Assume that there are $n$ buses. Fix $N>n$ and $T>0$. At time $t=0$ the $i^{th}$ bus is at location $1-i$, $1\leq i\leq n$. The evolution of the bus process $L_t=(L^{(1)}_t,\cdots, L^{(n)}_t)$,  is given by $n$ independent rate 1 Poisson processes conditioned not to intersect for $t\in [0,T]$, and subject to the terminal condition $L_T= (N, N-1,...,N-n+1)$. Consider buses arriving at a fixed point $x\in\{1,2,..., N-n+1\}$ with consecutive arrival times $0<t_1<\cdots <t_k <T$ (see Figure 1). 
\begin{figure}
 \centerline{\epsfig{file=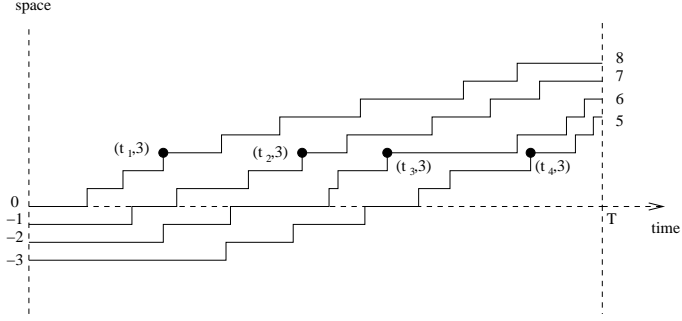, width=9cm}}
 \caption{Conditioned Poisson walks arriving at location 3}\label{figure1}
\end{figure}

The probability density for the arrival times $t_1,...,t_n$, $A_x(t_1,...,t_n)dt_1\cdots dt_n$, is given by  
\begin{equation}\label{conditionalprobability}
\frac{P\left[ \{(0,1-i)\} \rightarrow \{(t_i,x-1)\}; no\hspace{.1cm}intersection
\right] P\left[ \{(t_i,x)\}\rightarrow \{(T, N-(1-i))\}; no\hspace{.1cm}intersection
\right]dt_1 \cdots dt_n}{\mathbb{P} \left[ \{(0,1-i)\}\rightarrow \{(T, N-(1-i))\};  no\hspace{.1cm}intersection  \right]}.
\end{equation}
The fact that the numerator in \eqref{conditionalprobability} splits into a product of two factors is a consequence of the strong Markov property of $L_t$: in particular, functionals of the parts of the paths $\{ (t_i,x) \rightarrow (T, N-(1-i)), 1\leq i \leq n\}$ depend only on the increments of the process after some stopping times and hence, are independent of any functionals of the initial parts $\{ (0,1-i) \rightarrow (t_i, x-1), 1\leq i \leq n \}$ of the paths. Note also that the endpoints $(t_i,x-1)$, $1\leq i \leq n$, reflect the fact that the $i^{th}$ path jumps at $t_i$, and not before. Similar dynamical models with $T=\infty$ were considered in~\cite{KonigOConnellRoch}.

As we now show, each term in the ratio \eqref{conditionalprobability} has a determinantal form of Karlin-McGregor type~\cite{KarlinMcGregor}.
\begin{lemma}\label{KarlinMcGregor}
\begin{equation}
P\left[ \{(0,1-i)\} \rightarrow \{(t_i,x-1)\}; 
  no\hspace{.1cm}intersection \right] = det\left( e^{-t_j}
  \frac{t_j^{x+i-2}}{(x+i-2)!} \right)_{i,j=1}^n , \label{first} \\
\end{equation}
\begin{equation}
P\left[ \{(t_i,x)\}\rightarrow \{(T, N+(1-i))\};
  no\hspace{.1cm}intersection \right]= det \left( e^{-(T-t_j)}
  \frac{(T-t_j)^{N+(1-i)-x}}{(N+(1-i)-x)!} \right)_{i,j=1}^n, \label{second} \\
\end{equation}
\begin{equation}
\mathbb{P} \left[ \{(0,1-i)\}\rightarrow \{(T, N+(1-i))\};  no\hspace{.1cm}intersection  \right]=det\left(e^{-T} \frac{T^{N+i-j}}{(N+i-j)!} \right)_{i,j=1}^n.\label{third}
\end{equation}
\end{lemma}
\proof  \hspace{.2cm} Since the proofs of\eqref{first} and \eqref{second} equalities are similar, we only prove \eqref{first}. The proof is a variation on the standard Karlin-McGregor argument~\cite{KarlinMcGregor}. Fix $0<t_1< \cdots < t_n<T$ and define the following stopping
times:
\begin{eqnarray}
\tau_1&=&[inf\{s\in [0,t_1]: L_s^{(1)}=L_s^{(2)}\}]\wedge T \nonumber \\
\tau_2&=&[inf\{s\in [0,t_2]: L_s^{(2)}=L_s^{(3)}\}]\wedge T \nonumber \\
&&\cdots \nonumber \\
\tau_{n-1}&=&[inf\{s\in [0,t_{n-1}]: L_s^{(n-1)}=L_s^{(n)}\}]\wedge T. \nonumber
\end{eqnarray}
Observe that
\begin{eqnarray}
\det\left( e^{-t_j}\frac{t_j^{x+i-2}}{(x+i-2)!} \right)_{i,j=1}^n&=&
\sum_{\sigma\in S_n} sgn(\sigma) \prod_{i=1}^n
\mathbb{P}\left(L_{t_i}^{\sigma(i)}=x-1\right) \nonumber \\
&=&\sum_{\sigma\in S_n} sgn(\sigma) \mathbb{E} \left(\prod_{i=1}^n
\ind_{\{{L_{t_i}^{\sigma(i)}=x-1}\}}\right) \nonumber \\
&=&\sum_{\sigma\in S_n} sgn(\sigma) \mathbb{E}\left\{ \left[ \prod_{i=1}^n
  \ind_{\{{L_{t_i}^{\sigma(i)}=x-1}\}}\right] \left[
\ind_{\tau_1\wedge \cdots \wedge \tau_{n-1}=T} +
\sum_{j=1}^{n-1} \ind_{\tau_1\wedge \cdots \wedge
  \tau_{n-1}=\tau_j<T}\right]\right\} \nonumber \\
&=& \sum_{\sigma\in S_n} sgn(\sigma) \mathbb{E} \left[
  \ind_{\tau_1\wedge \cdots \wedge \tau_{n-1}=T} \prod_{i=1}^n
  \ind_{\{{L_{t_i}^{\sigma(i)}=x-1}\}}\right]\nonumber \\ 
&& + \sum_{j=1}^{n-1}\sum_{\sigma\in S_n} sgn(\sigma) \mathbb{E}  \left[ (\ind_{\tau_1\wedge \cdots \wedge
  \tau_{n-1}=\tau_j<T})\prod_{i=1}^n
  \ind_{\{{L_{t_i}^{\sigma(i)}=x-1}\}}\right] \nonumber \\
&=& \mathbb{P}\left[ \{(0,1-i)\} \rightarrow \{(t_i,x-1)\};
  no\hspace{.1cm}intersection \right] \nonumber \\ 
&& + \sum_{j=1}^{n-1}\sum_{\sigma\in S_n} sgn(\sigma) \mathbb{E}  \left[ (\ind_{\tau_1\wedge \cdots \wedge
  \tau_{n-1}=\tau_j<T})\prod_{i=1}^n
  \ind_{\{{L_{t_i}^{\sigma(i)}=x-1}\}}\right],  \label{determinant}
\end{eqnarray}
since the only nonzero term in the first sum is the term corresponding
to the identity permutation. We need only show that the second
summation vanishes. Let $\rho_j \in S_n$ be the transposition $(j,
j+1)$. By the strong Markov
property for $L_t$, 
\begin{equation}\nonumber
sgn(\sigma) \mathbb{E}\left[ (\ind_{\tau_1\wedge \cdots \wedge
  \tau_{n-1}=\tau_j<T})\prod_{i=1}^n
\ind_{\{{L_{t_i}^{\sigma(i)}=x-1}\}}\right] +sgn(\rho_j \sigma)
\mathbb{E}\left[ (\ind_{\tau_1\wedge \cdots \wedge
  \tau_{n-1}=\tau_j<T})\prod_{i=1}^n
\ind_{\{{L_{t_i}^{\rho_j \sigma(i)}=x-1}\}}\right] =0.
\end{equation}
Since for each $j$ the action of the given transposition $\rho_j$ is an involution on $S_n$,
\begin{equation}
\sum_{\sigma \in S_n} sgn(\sigma) \mathbb{E}\left[ (\ind_{\tau_1\wedge \cdots \wedge
  \tau_{n-1}=\tau_j<T})\prod_{i=1}^n
\ind_{\{{L_{t_i}^{\sigma(i)}=x-1}\}}\right] =0. \nonumber
\end{equation}
This immediately implies that the second sum in equation
\eqref{determinant} vanishes leaving only the desired
probability.

The proof of the third equality \eqref{third} is just the standard Karlin-McGregor argument~\cite{KarlinMcGregor}: In the proof of \eqref{first}, simply replace the $\tau_i$'s by
\begin{equation}
\tilde \tau_i=T\wedge \inf\{s\in[0,T]: L_s^i=L_s^{i+1}\} \text{,   for } i=1,...,n-1,  
\end{equation}
$x-1$ by $N+(1-i)$, and set $t_i=T$ for $1\leq i \leq n$. This completes the proof of the Lemma \ref{conditionalprobability}. \qed

The determinants in \eqref{first} and \eqref{second} are multiples of Vandermonde determinants. Indeed, simple computations show that
\begin{equation}\label{fourth}
\det\left( e^{-t_j} \frac{t_j^{x+i-2}}{(x+i-2)!} \right)_{i,j=1}^n = \frac{e^{-\sum_{j=1}^n t_j}}{\prod_{i=1}^n (x+i-2)!} \prod_{j=1}^n t_j^{x-1} \prod_{1\leq i < j \leq n}(t_j -t_i) 
\end{equation}
\begin{equation}\label{fifth}
\det \left( e^{-(T-t_j)} \frac{(T-t_j)^{N+(1-i)-x}}{(N+(1-i)-x)!} \right)_{i,j=1}^n= \frac{e^{-\sum_{j=1}^n (T-t_j)}}{\prod_{i=1}^{n} (N+1-x-i)!} \prod_{j=1}^n (T-t_j)^{N-x+1-n} \prod_{1\leq i < j \leq n}(t_j -t_i) 
\end{equation}
Inserting \eqref{third}, \eqref{fourth}, and \eqref{fifth} into \eqref{conditionalprobability} and setting $\tilde y_j=\frac{t_j}{T}$, one obtains
\begin{eqnarray}\label{FirstJacobiDensity}
A(t_1,...,t_n)dt_1 \cdots dt_n &=& \tilde A(\tilde y_1,..., \tilde y_n) d\tilde y_1 \cdots d\tilde y_n \nonumber \\
&=& C_{N,n,x} \prod_{i<j} (\tilde y_j -\tilde y_i)^2 \prod_{j=1}^n \tilde y_j^{x-1} (1-\tilde y_j)^{N-x-(n-1)}d\tilde y_1 \cdots d\tilde y_n,
\end{eqnarray} 
where $C_{N,n,x}=\left(\prod_{i=1}^n \frac{1}{(x+i-2)! (N-x+(1-i))!}\right)\left(\det\left(\frac{1}{(N+i-j)!}\right)_{i,j=1}^n\right)^{-1}$. The simple variable change $y_j= 2\tilde y_j -1$ induces the standard Jacobi weights on
$[-1,1]$:
\begin{equation} \label{StandardJacobi}
\tilde{C}_{N,n,x} \prod_{i<j} (y_j -y_i)^2 \prod_{j=1}^n
  (1+y_j)^{x-1} (1-y_j)^{N-x-(n-1)}dy_1 \cdots dy_n, 
\end{equation}
where $\tilde{C}_{N,n,x}=2^{-nN} C_{N,n,x}$. It is at this point that contact is made with RMT: If we view the $y_j$'s as eigenvalues of an $n\times n$ Hermitian matrix, then \eqref{StandardJacobi} is precisely the joint probability density for the eigenvalues of matrices in the so-called \underline{Jacobi Ensemble}. More precisely, for $n\times n$ Hermitian matrices $\{M=M^*\}$, set
\begin{equation}
w(M)=\det(1+M)^{x-1} \det(1-M)^{N-x-n+1} \det(\chi_{(-1,1)}(M)).
\end{equation}
Note that $w(M)$ is invariant under unitary conjugation of $M$, $M\rightarrow UMU^*$, for all unitary $U$. Then $P_n(M)dM\equiv \frac{1}{Z_n} w(M)dM$ defines the \underline{Jacobi Unitary Ensemble}, where $dM=\prod_{i=1}^n dM_{ii} \prod_{i<j} dM_{ij}^R \prod_{i<j} dM_{ij}^I$ and $Z_n$ is a normalization constant. A standard calculation (e.g. Mehta~\cite{Mehta}) then shows that  the distribution function for the eigenvalues, $\{\lambda_i\}$, of matrices in the ensemble is given by \eqref{StandardJacobi} with the identification $y_i=\lambda_i$. In summary, we see from \eqref{StandardJacobi} that for finite $N,n$ the arrival times for the buses at a location $x$ are distributed like the eigenvalues of a random matrix from the Jacobi Ensemble.  By techniques which are now  routine in random matrix theory (see~\cite{Mehta}),  \eqref{StandardJacobi} is amenable to asymptotic analysis in certain double scaling limits.  These limits (see Section \ref{Section3}) correspond to GUE: The fact that Jacobi Unitary Ensemble $\rightarrow$ GUE in the double scaling limit is a particular example of the well known phenomenon of \underline{universality} in RMT (see~\cite{DKMVZ1,DKMVZ2,KuVan}).

In the previous calculations we considered the arrival times $t_1,...,t_n$ for the buses at a fixed location, $x$. It is also of interest to consider the locations $x_1,...,x_n$ at a fixed recording time $t$. It turns out that the distribution of the $x_i$'s is again given by a random matrix type ensemble, the \underline{Krawtchouk Unitary Ensemble}.  In~\cite{KonigOConnellRoch}, K\"onig, O'Connell, and Roch discovered a similar relation between the Krawtchouk Unitary Ensemble and Poisson random walks conditioned not to intersect for all positive times. In analogy to \eqref{conditionalprobability}, we have that the probability $\mathcal{L}_t(x_1,...,x_n)$ for the positions $x_1>\cdots > x_n$ of the buses at time $t\in (0,T)$  is given by
\begin{equation}\label{secondratio}
\frac{\mathbb{P}\left( \{(0,1-i)\rightarrow (t,x_i)\}; \text{ no intersection} \right) \mathbb{P} \left( \{(t,x_i)\rightarrow (T, N+1-i) \}; \text{ no intersection} \right)}{\mathbb{P}\left(  \{(0,1-i)\rightarrow (T, N+1-i)\}; \text{ no intersection}     \right)},
\end{equation}
and by a Karlin-McGregor argument, similar to and even simpler than that in Lemma \ref{KarlinMcGregor}, we obtain the following:
\begin{lemma}
\begin{equation}\label{2.1}
\mathbb{P}\left( \{(0,1-i)\rightarrow (t,x_i)\}; \text{ no intersection} \right)=\det\left(e^{-t} \frac{t^{x_j+i-1}}{(x_j+i-1)!}\right)_{i,j=1}^n,
\end{equation}
\begin{equation}\label{2.2}
\mathbb{P} \left( \{(t,x_i)\rightarrow (T, N+1-i) \}; \text{ no intersection} \right)=\det\left(e^{(T-t)} \frac{(T-t)^{N+1-i-x_j}}{(N+1-i-x_j)!}\right)_{i,j=1}^n,
\end{equation}
\begin{equation}\label{2.3}
\mathbb{P}\left(  \{(0,1-i)\rightarrow (T, N+1-i)\}; \text{ no intersection}     \right)=\det\left(e^{-T} \frac{T^{N+i-j}}{(N+i-j)!}\right)_{i,j=1}^n.
\end{equation}
\end{lemma}
Inserting \eqref{2.1}, \eqref{2.2}, and \eqref{2.3} into \eqref{secondratio} we are led to the following formula for $\mathcal{L}_t(x_1,...,x_n)$~\cite{Johansson}:
\begin{eqnarray} \label{Lt}
\mathcal{L}_t(x_1,...,x_n)=&&[(N+n-1)!]^{-n}  \left[ \prod_{i=1}^n \frac{(N+i-1)!}{(i-1)!}\right]  (p-p^2)^{-\frac{n(n-1)}{2}} \prod_{i<j} (y_j - y_i)^2 \nonumber \\
&& \times \prod_{j=1}^n \frac{(N+n-1)!}{y_j! (N+n-1-y_j)!}  p^{y_j} (1-p)^{N+n-1-y_j},
\end{eqnarray}
where $y_j=x_j+n-1$, $p=\frac{t}{T}$, and $N+n-1\geq y_1>\cdots y_n\geq 0$. We have used the elementary identity 
\begin{equation}
\left[\det\left(\frac1{(N+i-j)!}\right)_{i,j=1}^n\right]^{-1} =\prod_{i=1}^n \frac{(N+i-1)!}{(i-1)!} 
\end{equation}
to derive \eqref{Lt}. As indicated above, \eqref{Lt} corresponds to the Krawtchouk
Unitary Ensemble.

An essential aspect in the analysis of the random particle ensembles \eqref{StandardJacobi} and \eqref{Lt} is that they are \underline{determinantal}, i.e. for each $r=1,2,...$, the $r^{th}$ correlation function for the ensemble can be expressed in determinantal form $\det(K(x_i,x_j))_{i,j=1}^r$ for some appropriate correlation kernel $K(x,y)$ (see for example, \cite{Mehta, Soshnikov, BorodinRains,BKPV}). It turns out that the fixed time distribution \eqref{Lt}, in particular, can be extended to a dynamical random particle system which is also determinantal. Indeed, fix $k$ times $0<t_1<\cdots < t_k<T$  and let $x^{(j)} = \{x_1^{(j)}, ..., x_n^{(j)} \}$, $1\leq j \leq k$, denote the positions of the buses at times $t_j$, $1\leq j\leq k$. Let $p_j=\frac{t_j}{T}$, $1\leq  j \leq k$. Then, using the Markov property of the bus system as above, we arrive at the following distribution for $x^{(1)},..., x^{(k)}$:
\begin{eqnarray}\label{multitime}
\text{Prob}[x^{(1)},...,x^{(k)}] = && const \cdot \prod_{i<j} (x_i^{(1)} - x_j^{(1)}) \prod_{i=1}^n \frac{p_1^{x^{(1)}_i}}{x^{(1)}_1 !} \nonumber \\ 
&& \times \det(V_{p_2 - p_1}(x_i^{(1)}, x^{(2)}_j)) \det(V_{p_3 - p_2}(x_i^{(2)}, x^{(3)}_j)) \cdots 
\det(V_{p_k - p_{k-1}}(x_i^{(k-1)}, x^{(k)}_j)) \nonumber \\
&& \times \prod_{i<j} (x^{(k)}_i - x^{(k)}_j ) \prod_{i=1}^n \frac{(1-p_k)^{N+n-1-x_j^{(k)}}}{(N+n-1-x_j^{(k)})!},
\end{eqnarray}
where  $V_p(x,y) = \frac{p^{y-x}}{(y-x)!}$ if $y\geq x$ and $V_p(x,y)=0$ otherwise. Using a variant of the Eynard-Mehta theorem~\cite{EynardMehta} as described in~\cite{BorodinOlshanski}, it follows directly that \eqref{multitime} is determinantal with an appropriate  correlation kernel $K_{t_i,t_j}(x,y)$. Further analysis yields the following block integral representation for the kernel: for $x\in x^{(i)}, y\in x^{(j)}$,
\begin{equation}
K_{t_i,t_j}(x,y)=\frac1{(2\pi i)^2} \oint \oint \frac{(p_j - (1-p_j)t)^y (1+t)^{N+n-1-y} t^{-n} }{(p_i-(1-p_i)s)^{x+1} (1+s)^{N+n-x} s^{-n}} \frac{dsdt}{t-s},
\end{equation}
 where the integration contours are chosen as follows: $s$ runs along a simple positively oriented contour which goes around $\frac{p_i}{1-p_i}$ and does not contain $-1$; $t$ runs along a simple positively oriented contour which goes around $0$; for $p_i \geq p_j$ the $s$-contour contains the $t$-contour and for $p_i <p_j$ the $s$-contour lies inside the $t$-contour. 
 
 In the following section we will analyze \eqref{StandardJacobi} and \eqref{Lt} asymptotically as $N,n,x,t,T\rightarrow \infty$. We plan to present the derivation and asymptotic analysis of  \eqref{multitime} in a later publication.

\section{The Double Scaling Limits} \label{Section3}

We will first describe the arrival statistics of the
buses \eqref{StandardJacobi} in the limit $N,n,x\rightarrow \infty$ where $\frac{n}{N}
\rightarrow \nu$, $\frac{x-1}{N} \rightarrow \eta$, and $0<\nu,\eta<1$. This scaling
limit is the natural one from both the physical and mathematical point
of view: the number of buses, $n$, should be proportional to the length
of the bus route, $N+n$, and the arrival point, $x$, should not be
too close to the beginning or the  end of the route. The main result is that in the ``unfolded'' scale (see below), the Christoffel-Darboux kernel (equivalently, the two point correlation function) $K_{N,n,x}$ below, converges as $N\rightarrow \infty$, $\frac{n}{N} \rightarrow \nu$, $\frac{x}{N} \rightarrow \eta$ to the so-called sine kernel (see below), which is \underline{universal} in random matrix theory. The proof of the convergence $K_{N,n,x} \rightarrow \mathcal{K}_\infty$ for the Jacobi ensemble was given first by Nagao and Wadati~\cite{NagaoWadati}.

In RMT the analysis of the spacing distribution proceeds via the analysis of the gap distribution (see, for example, ~\cite{Mehta, Deift}): 
\begin{equation} \label{gapprobability}
\mathbb{P}_{N,n}\left((c,d)\right)=\text{Prob}(y_j\notin(c,d) \text{, i=1,...,n}),
\end{equation}
There are three steps involved in the analysis of $\mathbb{P}_{N,n}\left((c,d)\right)$: (a)The first step is to express $\mathbb{P}_{N,n}((c,d))$ in Fredholm determinantal form, 
\begin{equation}
\mathbb{P}_{N,n}((c,d))=\det(I-K)_{L^2[c,d]},
\end{equation} 
for some kernel operator $K$ expressed in terms of certain appropriate orthogonal polynomials (see, in particular, Tracy and Widom~\cite{TracyWidom}); (b) the second step of the
analysis is to determine the appropriate scalings for $c$ and $d$;
this is accomplished by evaluating  the so-called equilibrium measure for
the associated logarithmic potential theoretic problem~\cite{SaffTotik, Deift}; (c) the third
step involves a detailed asymptotic analysis of the appropriate associated
orthogonal polynomials in the double scaling limit.

(a) For $j=1,2,...$ , let $p^{N,n,x}_j(y) = \gamma_j^{N,n,x} y^j +...$, $\gamma_j^{N,n,x}>0$ denote the Jacobi polynomials obtained by orthonormalizing $\{1,y,y^2,...\}$ with respect to the weight $w_{N,n,x}(y)= (1-y)^{N-(n+x-1)}(1+y)^{x-1}\chi_{[-1,1]}(y)$. Set  
$\phi_j^{N,n,x} (y)=p_j^{N,n,x}(y) w_{N,n,x}^\frac12 (y)$, for $j=1,2,...$. Then,  the probability that there is no bus
arrival in the interval of time $[a_{N,n,x}, b_{N,n,x}]$, $\mathbb{P}_{N,n}((a_{N,n,x}, b_{N,n,x}))$, is given by  
\begin{equation}
det(I-K_{N,n,x})_{L^2([a_{N,n,x}, b_{N,n,x}])},
\end{equation}
where 
\begin{equation}
K_{N,n,x}(z,z')=l_n\frac{\phi_n^{N,n,x} (z)\phi_{n-1}^{N,n,x}
  (z')-\phi_n^{N,n,x} (z')\phi_{n-1}^{N,n,x} (z)}{z-z'},
\end{equation}
and $l_n$ is given by the formula $\frac{k_n}{k_{n+1}} \sqrt{\frac{h_{n+1}}{h_n}}$ where $h_n$ and $k_n$ are constants arising in the three term recurrence relation for general orthogonal polynomials; for the Jacobi polynomials,  $h_n$ and $k_n$ are given, respectively,  by formulas 22.2.1 and 22.3.1 of Abramowitz and Stegun~\cite{AbromowitzStegun} in the case $\alpha=N-(n+x-1)$ and $\beta=x-1$.

(b) By standard methods (see Saff and Totik~\cite{SaffTotik} and Deift~\cite{Deift}), the equilibrium measure can be determined explicitly for the problem at hand. It turns out that the support of the equilibrium measure is an interval $[a,b]\subset [-1,1]$, and for 
fixed $\nu, \eta>0$  satisfying  $\nu+\eta<1$, the 
measure takes the form:
\begin{equation}\label{EquilibriumMeasure}
\psi(x)dx=\frac{\eta \sqrt{(x-a)(b-x)}}{\pi \sqrt{(1+a)(1+b)} (1-x^2)}dx,
\end{equation}
where $a$ and $b$ satisfy the relations:
\begin{equation} \nonumber
\frac{\eta}{\sqrt{(1+a)(1+b)}} = -\frac{(\eta +\nu
  -1)}{\sqrt{(1-a)(1-b)}}, \hspace{.8cm}
(1+\nu)=\frac{\eta}{\sqrt{(1+a)(1+b)}}- \frac{(\eta +\nu -1)}{\sqrt{(1-a)(1-b)}},
\end{equation}
and the radicals represent the positive square root. In the symmetric
case where $\eta=\frac{1-\nu}{2}$,
\begin{equation} \nonumber
b=-a=\left(1-\left(\frac{1-\nu}{1+\nu}\right)^2 \right)^\frac{1}{2}.
\end{equation}

(c) Let $\tau_0 \in (a,b)$. A standard calculation in RMT show that the expected number of particles per unit interval in a neighborhood of $\tau_0  \approx n\psi(\tau_0) $. Thus, changing
scales $t \mapsto s:= n\psi(\tau_0)t$, we see that the expected number of
particles per unit s-interval is 1: this process of rescaling is known as
"unfolding" the data. The probability that there are no bus arrivals in $[\tau_0-\frac{s}{\psi(\tau_0)
  n},\tau_0+\frac{s}{\psi(\tau_0) n}]$ becomes
\begin{equation}
det\left(I-\mathcal{K}_n\right)_{L^2(-s,s)},
\end{equation}
where 
\begin{equation}
\mathcal{K}_n(\tau_0,\xi, \rho)= l_n \frac{\phi_n^{N,n,x} (\tau_0 +
  \frac{\xi}{\psi(\tau_0)n})\phi_{n-1}^{N,n,x}
  (\tau_0 + \frac{\rho}{\psi(\tau_0)n})-\phi_n^{N,n,x} (\tau_0 +
  \frac{\rho}{\psi(\tau_0)n})\phi_{n-1}^{N,n,x} (\tau_0 +
  \frac{\xi}{\psi(\tau_0)n})}{\xi-\rho}.
\end{equation}
Analyzing the asymptotics of the Jacobi polynomial as $N,n,x\rightarrow \infty$ as above, we finally see that 
\begin{equation}
det\left(I-\mathcal{K}_n\right)_{L^2(-s,s)} \rightarrow det\left(I-\mathcal{K}_\infty\right)_{L^2(-s,s)},
\end{equation}
where $\mathcal{K}_\infty(\xi, \rho) =
\frac{sin\pi(\xi-\rho)}{\pi(\xi-\rho)}$, the so-called sine kernel.  To compute the limiting conditional
probability that given a bus arrival at time $s$, the next bus arrives
at time $t$ we simply compute (see, e.g.~\cite{TracyWidom})
\begin{equation}\label{conditional}
\lim_{N\rightarrow \infty} \frac{1}{\psi(s)} \frac{\partial^2}{\partial s \partial t}det\left(I-\mathcal{K}_N \right)_{L^2(s,t)}= \frac{1}{\psi(s)} \frac{\partial^2}{\partial s \partial t}det\left(I-\mathcal{K}_\infty\right)_{L^2(s,t)}.
\end{equation}
Formula \eqref{conditional}, the Gaudin distribution,  gives the exact spacing distribution of RMT; as noted in the Introduction, this distribution in known to be well approximated by the Wigner surmise used by \v{K}rbalek and \v{S}eba~\cite{KS} to analyze their observations of the bus arrival times. 

In addition, using the fact that $K_{N,n,x} \rightarrow K_\infty$ in the unfolded scale, one easily sees that the number variance, $H_{N,n,x}(\cdot)$, say, converges i.e., 
\begin{equation}
H_\infty (s) := \lim_{N\rightarrow \infty, \frac{n}{N} \rightarrow \nu, \frac{x-1}{N} \rightarrow \eta}   H_{N,n,x}\left(\frac{s}{n\psi(0)}\right)
\end{equation}
exists. Moveover, $H_\infty (s)$ is precisely the number variance for GUE, as given, for example, in formula (16.1.3) in~\cite{Mehta}. As $s$ becomes large, $H_\infty (s) \sim \frac1{(\pi)^2} (\log(2\pi s) +\gamma +1)$ where  $\gamma$ is Euler's constant, as noted in formula (5) of~\cite{KS}.

The analysis of the fixed time bus locations \eqref{Lt} as $N,T,t\rightarrow \infty$,  $\frac{t}{T}\rightarrow C_1$, $\frac{n}{N} \rightarrow C_2$, and $\frac{T}{N}\rightarrow C_3$ is very similar and follows the same general procedure (a), (b), (c) above. We note that the Krawtchouk ensemble was first analyzed in the same scaling limit by Johansson in his analysis of the Aztec diamond~\cite{Johansson}. Again, one finds that in the limit the statistics of the bus locations is governed by the sine kernel.

\section{The Circular Bus Route}\label{Section4}

In the previous bus model, we considered the buses moving from one terminus at $t=0$ to a second terminus at time $t=T$.  In this final section, we consider a bus model on a circular route. K\"onig, O'Connell, and Roch(see again~\cite{KonigOConnellRoch}) investigated a related model on the discrete circle. Let
$\mathbb{Z}_M$ be the discrete circle with nodes labeled by
$\{0,1,2,...,M-1\}$. We will consider the case of $k<M$ buses
traveling along $\mathbb{Z}_M$. The buses  start at positions
$0\leq N_1(0) <\cdots <N_k(0)\leq M-1$ and evolve as independent Poisson processes
conditioned not to intersect. Let $\theta, \tilde \theta\in
\mathbb{Z}_M$ and note that the transition probability on $\mathbb{Z}_M$ for a single
rate 1 Poisson process to travel from state $\theta$ to $\tilde
\theta$ is  
\begin{equation}
p_t(\theta, \tilde \theta)= e^{-t}\sum_{l=-\infty}^\infty
\frac{t^{\tilde \theta -\theta +lM}}{(\tilde \theta -\theta +lM)!}.
\end{equation}
Let $(\theta_1,...,\theta_k)$ be a $k$-tuple of distinct elements of
$\mathbb{Z}_M$ such that $0\leq \theta_1<...<\theta_k \leq M-1$. Let $(\tilde
\theta_1,...,\tilde \theta_k)$ be another $k$-tuple of distinct elements of
$\mathbb{Z}_M$ such that there exists a cyclic permutation $\sigma\in
S_k$ for which $(\tilde \theta_{\sigma(1)}, ... , \tilde \theta_{\sigma(k)})$
satisfies $0\leq \tilde \theta_{\sigma(1)}<...<\tilde \theta_{\sigma(k)} \leq M-1$. Given
that the buses begin at positions $\theta_1<...<\theta_k$ at time
$0$, the probability that they are at locations $\tilde
\theta_1,...,\tilde \theta_k$ at time t, and have not intersected in
the mean time, is given by the determinantal expression
\begin{equation}
\det\left(p_t(\theta_i, \tilde \theta_j)\right)_{i,j=1}^k. \label{CircleKM}
\end{equation}
We prove this expression by adapting an alternate proof of the Karlin-McGregor formula from the line to the circle as follows.

Let $\mathbb{Z}^k_{M\geq}$ be defined as the subset of $\mathbb{Z}^k /M\mathbb{Z}^k$ for which there exists a cyclic permutation $\sigma \in S_d$  such that $0\leq z_{\sigma(1)} \leq z_{\sigma(2)} \leq \cdots \leq z_{\sigma(k)}\leq M-1$. For $f:\mathbb{Z}^k_{M\geq} \rightarrow \mathbb{R}$, define the operator
\begin{equation}
Lf(z)=\sum_{i=1}^k [f(z_1,...,z_i +1,...,z_k) - f(z)] \text{ , for } z\in \mathbb{Z}^k_{M\geq} - \partial  \mathbb{Z}^k_{M\geq} .
\end{equation}

As we will see, in order to prove the circular Karlin-McGregor formula \eqref{CircleKM}, it is sufficient to show that the function $g_{\tilde \theta}:[0,t]\times \mathbb{Z}^k_{M\geq}\rightarrow \mathbb{R}$ defined by $g_{\tilde\theta}(t,\theta)= \det(p_{t-s}(\theta_i,\tilde \theta_j))$ is the unique solution to the equation 
\begin{eqnarray}\label{backwardequation}
&&(\partial_s + L)g_{\tilde \theta}=0, \\
&&g_{\tilde \theta}|_{[0,t]\times \partial \mathbb{Z}^k_{M\geq}}=0, \\
&&g_{\tilde \theta}(t,\theta)=\delta_{\tilde \theta}(\theta).\label{terminalcondition}
\end{eqnarray}
 One checks easily that  $g_{\tilde \theta}$ is a solution of \eqref{backwardequation}. In order to check that the solution is unique, one needs to prove a maximum principle for this equation, but this is easily done by mimicking the standard proof for parabolic equations.
 
We now show that $g_{\tilde \theta}(0,\theta)$ gives the desired probability. Let $N(s)=(N_1(s),...,N_k(s))$ be the Poisson process on $\mathbb{Z}^k_{M\geq}$ with initial condition $N(0)=(\theta_1,...,\theta_k)$.  Let $\tau=\inf\{s\in [0,t]: N(s)\in \partial \mathbb{Z}^k_{M\geq} \}$. By Ito's formula, 
\begin{equation}
g_{\tilde \theta}(s,N(s\wedge \tau) - g_{\tilde \theta}(0,\theta) - \int_0^{s \wedge \tau} (\partial_r + L)g_{\tilde \theta}(r, N(r))dr
\end{equation}
is a martingale. Taking expectations, we obtain
\begin{equation}\label{circnonint}
g_{\tilde \theta}(0,\theta)=\mathbb{E}^\theta g(s, N({s\wedge \tau}) \text{ for } s \in [0,t].
\end{equation}
Letting $s=t$, \eqref{circnonint} becomes $g_{\tilde \theta}(0,\theta)=\mathbb{E}^\theta g(t, N(t\wedge \tau)$. Using the definition of $\tau$ and the terminal condition \eqref{terminalcondition}, we see that 
\begin{equation} \label{nonintfinal}
g_{\tilde \theta}(0,\theta)=\mathbb{P}^\theta (t<\tau \text{ and } N_t=\tilde \theta),
\end{equation}
which proves \eqref{CircleKM}.  Thus, for buses starting at locations $0\leq \theta_1 < \cdots < \theta_k \leq M-1$ at time $t=0$, we have 
\begin{eqnarray}
&&\text{Prob}(\text{busses are at locations} 0\leq \tilde \theta_1< \cdots \tilde \theta_k \leq M-1\text{ at time t }; \text{ no intersection for all } 0\leq s \leq t) \nonumber \\
&& = \mathbb{P} (t<\tau \text{ and } N(t)=\tilde \theta) = g_{\tilde \theta}(0, \theta) = \det(p_t(\theta_i, \tilde \theta_j))_{i,j=1}^k , \nonumber
\end{eqnarray}
as desired.

This formula enable us to obtain, in particular, a Karlin-McGregor type formula for the solution of the following natural problem for the buses on a circular bus route. Imagine that the initial locations of the buses are at $\theta_i=i-1$ for $i=1,...,k$. Suppose that the buses return to these locations at some fixed time $T$ later. For any time $0<t<T$, it immediately follows from the formula above and the Markov property that the distribution, $Q_t(\tilde \theta_1,...,\tilde \theta_k)$, of the locations of the buses conditioned on arriving at $\theta_1,....,\theta_k$ at time $T$, is given by the formula:
\begin{equation}
Q_t(\tilde \theta_1,...,\tilde \theta_k)=\frac{\det(p_t(\theta_i, \tilde \theta_j)\det(p_{T-t}(\tilde \theta_i, \theta_j))}{\det(p_T(\theta_i, \theta_j)},
\end{equation}
into which we may now substitute  \eqref{CircleKM}. Observe that the number of rotations about the circle is not fixed. We plan to investigate the asymptotic behavior of this model in a future paper.


\begin{thebibliography}{10}

\bibitem{AbromowitzStegun}
M.~Abramowitz and I.~Stegun, editors.
\newblock {\em Handbook of mathematical functions with formulas, graphs, and
  mathematical tables}.
\newblock Dover Publications Inc., New York, 1992.
\newblock Reprint of the 1972 edition.

\bibitem{BKPV}
J.~Ben~Hough, M.~Krishnapur, Y.~Peres, and B.~Virag.
\newblock Determinantal processes and independence.
\newblock {\em arxiv.org/abs/math.PR/0503110}.

\bibitem{BorodinOlshanski}
A.~Borodin and G.~Olshanski.
\newblock Markov processes on partitions.
\newblock {\em arxiv.org/abs/math/0409075, to appear in Prob. Th. Rel. Fields}.

\bibitem{BorodinRains}
A.~Borodin and E.~Rains.
\newblock Eynard-mehta theorem, schur process, and their pfaffian analogs.
\newblock {\em arxiv.org/abs/math-ph/0409059, to appear in J. Stat. Phys.}


\bibitem{DKMVZ2}
P.~Deift, T.~Kriecherbauer, K.~T-R McLaughlin, S.~Venakides, and X.~Zhou.
\newblock Strong asymptotics of orthogonal polynomials with respect to
  exponential weights.
\newblock {\em Comm. Pure Appl. Math.}, 52(12):1491--1552, 1999.

\bibitem{DKMVZ1}
P.~Deift, T.~Kriecherbauer, K.~T.-R. McLaughlin, S.~Venakides, and X.~Zhou.
\newblock Uniform asymptotics for polynomials orthogonal with respect to
  varying exponential weights and applications to universality questions in
  random matrix theory.
\newblock {\em Comm. Pure Appl. Math.}, 52(11):1335--1425, 1999.

\bibitem{Deift}
P.~A. Deift.
\newblock {\em Orthogonal polynomials and random matrices: a
  {R}iemann-{H}ilbert approach}, volume~3 of {\em Courant Lecture Notes in
  Mathematics}.
\newblock New York University Courant Institute of Mathematical Sciences, New
  York, 1999.

\bibitem{EynardMehta}
B.~Eynard and M.~Mehta.
\newblock Matrices coupled in a chain. {I}. {E}igenvalue correlations.
\newblock {\em J. Phys. A}, 31(19):4449--4456, 1998.

\bibitem{Johansson}
K.~Johansson.
\newblock Non-intersecting paths, random tilings and random matrices.
\newblock {\em Probab. Theory Related Fields}, 123(2):225--280, 2002.

\bibitem{KarlinMcGregor}
S.~Karlin and J.~McGregor.
\newblock Coincidence probabilities.
\newblock {\em Pacific J. Math.}, 9:1141--1164, 1959.

\bibitem{KonigOConnell}
W.~K{\"o}nig and N.~O'Connell.
\newblock Eigenvalues of the {L}aguerre process as non-colliding squared
  {B}essel processes.
\newblock {\em Electron. Comm. Probab.}, 6:107--114 (electronic), 2001.

\bibitem{KonigOConnellRoch}
W.~K{\"o}nig, N.~O'Connell, and S.~Roch.
\newblock Non-colliding random walks, tandem queues, and discrete orthogonal
  polynomial ensembles.
\newblock {\em Electron. J. Probab.}, 7:no. 5, 24 pp. (electronic), 2002.


\bibitem{KS}
M.~\v{K}rbalek and P.~\v{S}eba.
\newblock Statistical properties of the city transport in Cuernevaca (Mexico)
  and random matrix theory.
\newblock {\em J. Phys. A: Math. Gen.}, 33:229--234, 2000.

\bibitem{KuVan}
A.~B.~J. Kuijlaars and M.~Vanlessen.
\newblock Universality for eigenvalue correlations from the modified {J}acobi
  unitary ensemble.
\newblock {\em Int. Math. Res. Not.}, (30):1575--1600, 2002.

\bibitem{Mehta}
M.~Mehta.
\newblock {\em Random matrices}, volume 142 of {\em Pure and Applied
  Mathematics (Amsterdam)}.
\newblock Elsevier/Academic Press, Amsterdam, third edition, 2004.

\bibitem{NagaoWadati}
T.~Nagao and M.~Wadati.
\newblock Correlation functions of random matrix ensembles related to classical
  orthogonal polynomials.
\newblock {\em J. Phys. Soc. Japan}, 60(10):3298--3322, 1991.

\bibitem{SaffTotik}
E.~Saff and V.~Totik.
\newblock {\em Logarithmic potentials with external fields}, volume 316 of {\em
  Grundlehren der Mathematischen Wissenschaften [Fundamental Principles of
  Mathematical Sciences]}.
\newblock Springer-Verlag, Berlin, 1997.
\newblock Appendix B by Thomas Bloom.

\bibitem{Soshnikov}
A.~Soshnikov.
\newblock Determinantal random point fields.
\newblock {\em Uspekhi Mat. Nauk}, 55(5(335)):107--160, 2000.

\bibitem{TracyWidom}
C.~Tracy and H.~Widom.
\newblock Correlation functions, cluster functions, and spacing distributions
  for random matrices.
\newblock {\em J. Statist. Phys.}, 92(5-6):809--835, 1998.


\end{thebibliography}

\end{document}